\documentclass[final,3p,times]{elsarticle}
\usepackage{amsmath,amssymb,bbm,bm,comment,enumitem,hyperref,mathrsfs,algorithm,algorithmic,soul,dsfont,xcolor}
\newtheorem{proposition}{Proposition}
\newtheorem{theorem}{Theorem}
\newtheorem{lemma}{Lemma}

\newtheorem{remark}{Remark}

\newboolean{showcomments}
\setboolean{showcomments}{true}
\DeclareRobustCommand{\red}[1]{\ifthenelse{\boolean{showcomments}}{\textcolor{black}{#1}}{}}
\DeclareRobustCommand{\blue}[1]{\ifthenelse{\boolean{showcomments}}{\textcolor{black}{#1}}{}}
\DeclareRobustCommand{\magenta}[1]{\ifthenelse{\boolean{showcomments}}{\textcolor{black}{#1}}{}}

\journal{Operations Research Letters}

\begin{document}

\begin{frontmatter}
	
\title{Some Asymptotic Properties of the Erlang-C Formula\\ in Many-Server Limiting Regimes\\ }
\author[label1]{Ragavendran Gopalakrishnan\fnref{fn1}}\ead{r.gopalakrishnan@queensu.ca}
\author[label2]{Yueyang Zhong\corref{cor1}}\ead{yzhong0@chicagobooth.edu}
\cortext[cor1]{Corresponding author}
\fntext[fn1]{This work is supported in part by funding from the Social Sciences and Humanities Research Council of Canada.}
\affiliation[label1]{organization={Smith School of Business at Queen's University},
    addressline={143 Union Street West},
    city={Kingston},
    postcode={ON K7L 3N6},
    country={Canada}}
\affiliation[label2]{organization={The University of Chicago Booth School of Business},
    addressline={5807 S Woodlawn Ave},
    city={Chicago},
    postcode={IL 60637},
    country={United States}}

\begin{abstract}
	This paper presents asymptotic properties of the Erlang-C formula in a spectrum of many-server limiting regimes. Specifically, we address an important gap in the literature regarding its limiting value in critically loaded regimes by studying extensions of the well-known square-root safety staffing rule used in the Quality-and-Efficiency-Driven (QED) regime.
\end{abstract}

		
\end{frontmatter}

\section{Introduction} \label{section:intro}

Multiserver systems are widely used to model situations in which customers may be served by one among multiple servers. Classical examples of such systems include call centers~\cite{erlang1917solution,whitt1999dynamic,gans2003telephone,brown2005statistical,aksin2007modern}, healthcare delivery~\cite{green2007providing,yom2010queues,armony2015patient}, and communication systems~\cite{anick1982stochastic,kelly1985stochastic}. The most basic multiserver queueing model is the $M$/$M$/$N$ queue (also known as the Erlang delay model), where customers arrive according to a Poisson process with rate $\lambda$ and are served by one of $N$ parallel servers for an exponentially distributed amount of time with mean $1/\mu$. A newly arriving customer that finds all servers busy joins a first-come-first-served queue and waits for their turn. Let $\rho = \lambda/\mu$ be the offered load.

A fundamental performance measure of any queueing model is the steady-state probability of delay, i.e., the steady-state probability that an incoming customer does not find an available server immediately upon entry and must therefore wait for service. For the $M$/$M$/$N$ queue, this quantity is given by the well-known Erlang-C formula~\citep[p.91]{Cooper72}:
\begin{align}\label{eq:ErlangC}
	C\left(N,\rho\right) = \frac{\rho^N}{N!} \left(\left(1 - \frac{\rho}{N}\right) \sum_{i=0}^{N-1}\frac{\rho^i}{i!} + \frac{\rho^N}{N!}\right)^{-1}.
\end{align}
The Erlang-C formula is heavily relied upon in a wide range of optimization problems in many-server queueing systems, such as optimal staffing problems; see, e.g.,~\cite{gans2003telephone,borst2004dimensioning} and references within. However, finding closed-form solutions for such problems is intractable due to the complexity of the Erlang-C formula. One approach to tackle this challenge is by developing  approximations for finite systems~\cite{harel1988sharp,kolesar1998insights,turpin2023}. However, more accurate approximations tend to be more analytically complicated, which can be a disadvantage of this method~\cite{turpin2023}. Alternatively, motivated by large-scale service systems, many-server limits (as $N$ and $\rho$ grow large while $\mu$ remains fixed) have been used to develop analytically simpler approximations~\cite{van2019economies}.

For the $M$/$M$/$N$ queue, the seminal paper~\cite{halfin1981heavy} shows that, when $N$ and $\rho$ grow unboundedly according to the relationship $N = \rho + \beta \sqrt{\rho}$ for some fixed $\beta > 0$, the steady-state probability of delay, namely, the Erlang-C formula~\eqref{eq:ErlangC}, converges to a value that is strictly between~$0$ and~$1$. This relationship is known as the square-root safety staffing rule. It achieves high system utilization, since $\rho/N$ grows closer to its critical value~$1$ (and therefore this staffing rule belongs to the class of critically loaded staffing rules), yet short customer waiting time on the order of $1/\sqrt{N}$~\cite{halfin1981heavy}. In other words, the square-root safety staffing rule achieves a balance between the dual goals of system efficiency and quality of service. Thus, large-scale systems under the square-root staffing rule are said to operate in the Quality-and-Efficiency-Driven (QED) many-server heavy-traffic limiting regime. Henceforth, we use the terms ``square-root safety staffing rule'' and ``QED regime'' interchangeably.

The QED regime is, in fact, the only limiting regime in which the steady-state probability of delay in the $M$/$M$/$N$ queue admits a non-degenerate limit (i.e., a limit that is neither $0$ nor $1$)~\cite{halfin1981heavy}. Furthermore, it is asymptotically optimal to operate the $M$/$M$/$N$ queue in the QED regime for large and heavily loaded systems, when choosing the optimal staffing level that minimizes a linear combination of staffing and waiting costs or when choosing the smallest staffing level subject to an upper bound on the waiting cost~\cite{borst2004dimensioning}.

%$\beta>0$ was required in the original QED regime (\cite{halfin1981heavy}) to guarantee stability when there is no abandonment. Our analysis, however, covers all values of $\beta$.
%For the $M$/$M$/$N$ queue, \cite{borst2004dimensioning} show that square-root safety staffing, first introduced in~\cite{erlang1948rational} and later formalized in~\cite{halfin1981heavy}, is asymptotically optimal for the aforementioned staffing problem, ...
%The QED regime is also widely regarded as one of the most suitable heavy-traffic limiting regimes in which to approximate the performance of variants of the basic $M$/$M$/$N$ system, such as many-server queueing systems with customer abandonment (\cite{garnett2002designing,zeltyn2005call}), call centers with a call-back option (\cite{armony2004contact,armony2004customer}), and call centers with multiple customer classes and agent skills (\cite{harrison2004dynamic,gurvich2008service,gurvich2009queue}).
%in that large-scale systems that operate in the QED regime dwarf the usual tradeoff between high system utilization and short waiting times. %balanace the service capacity and traffic demand so as to achieve a certain target performance standard or optimize a certain cost criterion.	
However, there exist several other well-motivated objectives in the $M$/$M$/$N$ queue and related variants for which the square-root safety staffing is not always asymptotically optimal~\cite{mandelbaum2009staffing,gopalakrishnan2016routing,zhan2019staffing,zhong2022asymptotically}, even among critically loaded staffing policies~\cite{nair2016provisioning,kim2018value}. 

In this paper, we address an important gap in the literature regarding the limiting value of the Erlang-C formula under critically loaded staffing rules, by studying extensions of the well-known square-root safety staffing rule, $N=\rho+\beta\sqrt{\rho}$. Specifically, we consider staffing rules of the form $N=\rho+f(\rho)$, where $f$ is a sublinear function that can be positive or negative. Compared to the square-root safety staffing rule, we make the following two inclusions concerning its second-order term: 
\begin{enumerate}
	\item [(i)] It can be any sublinear term.
	\item [(ii)] It can be negative, i.e., the system can be understaffed (resulting in insufficient capacity to meet the offered load). This is motivated by systems in which customers are lost because they are turned away upon arrival due to fully occupied waiting rooms and/or because they become impatient while in the system and leave before their service is completed (see~\cite{gans2003telephone} and the references therein). 
\end{enumerate}

In doing so, we identify more staffing rules (in addition to the square-root safety staffing rule) under which the Erlang-C formula admits a non-degenerate limit, making studying the limiting properties of key performance measures of these staffing rules more tractable; see Section~\ref{sec:implication}. This could potentially improve optimal system design by aiding the exploration of an expanded set of candidate staffing rules using many-server heavy-traffic approximations.

Our result unifies all the many-server limiting regimes, including the aforementioned QED regime, the Efficiency-Driven (ED) regime (where $\rho/N$ becomes larger than~$1$ in the limit, indicating an overloaded system, and the Quality-Driven (QD) regime (where $\rho/N$ becomes smaller than~$1$ in the limit, indicating an underloaded system).
\\
%We further illustrate the value and potential of our results by using them to (1) significantly simplify the proof of Theorem 1 in~\cite{nair2016provisioning}; and (2) formally extend the results in~\cite{kim2018value} to the many-server setting. \\

%Moreover, our results can significantly simplify the proof of Theorem~1 in~\cite{nair2016provisioning}, which determines the optimal capacity provisioning for a profit-maximizing firm. For the interested reader, we provide the simplified proof in~\ref{sec:appendix-new-proof-of-nair}.

\noindent\textbf{Notation.} We conclude this section by introducing some notations that will be used throughout the paper. We use the $o$, $\omega$, and $\Theta$ notations to denote the limiting behavior of functions. Formally, for any two real-valued functions $f(x)$ and $g(x)$ that take nonzero values for all sufficiently large $x$, we say that $f(x) \in o\left(g(x)\right)$ (equivalently, $g(x) \in \omega\left(f(x)\right)$) if $\lim_{x \to \infty} \frac{f(x)}{g(x)} = 0$, and $f(x) \in \Theta\left(g(x)\right)$ if $\lim_{x \to \infty} \left|\frac{f(x)}{g(x)}\right| \in (0,\infty)$. Moreover, the relation $f(x) \sim g(x)$ means $\lim_{x\to\infty} \frac{f(x)}{g(x)} = 1$.
Let $\phi(x) = \frac{1}{\sqrt{2\pi}} e^{-\frac{1}{2}x^2}$ and $\Phi^c(x) = \int_{x}^{\infty} \phi(t) dt$ denote the probability density function and the complementary cumulative distribution function of the standard normal distribution, respectively. Finally, define $\xi(x) := x\frac{\Phi^c(x)}{\phi(x)}$ for $x \in \mathbb{R}$.

\section{Main Result}\label{sec:main-result}
Consider a sequence of systems indexed by the arrival rate $\lambda$, and let $\lambda$ become large. Our convention, when we refer to any process or quantity associated with the system having arrival rate $\lambda$, is to superscript the appropriate symbol by $\lambda$. For example, we denote by $N^\lambda$ the number of servers in the system with arrival rate $\lambda$, and $\rho^{\lambda}=\lambda/\mu$ the corresponding offered load. We are interested in many-server limiting regimes obtained by letting the arrival rate $\lambda$ and the number of servers $N^{\lambda}$ grow unboundedly while the service rate $\mu$ remains fixed; nevertheless, our result continues to hold when $N^\lambda$ remains bounded as it grows with $\lambda$.

\begin{theorem}\label{thm:ErlC-ASYM}
	\begin{align*}
		\lim_{\lambda\to\infty} C\big(&N^{\lambda},\rho^{\lambda}\big) \\
		=&
		\begin{cases}
			\infty,&  0 < \rho^{\lambda} - N^{\lambda} \in \omega(\!\sqrt{\lambda}), \\
			\big(1 - \xi(z)\big)^{-1} \in (1,\infty), &  0 < \rho^{\lambda} - N^{\lambda} \in \Theta(\!\sqrt{\lambda}), \\
			1, &  |N^{\lambda} - \rho^{\lambda}| \in o(\!\sqrt{\lambda}), \\
			\big(1 - \xi(z)\big)^{-1} \in (0,1), & 0 < N^{\lambda} - \rho^{\lambda} \in \Theta(\!\sqrt{\lambda}), \\
			0, & 0 < N^{\lambda} - \rho^{\lambda} \in \omega(\!\sqrt{\lambda}),
		\end{cases}
	\end{align*}
	where $z = \lim_{\lambda\to\infty} \frac{\rho^{\lambda} - N^{\lambda}}{\sqrt{N^{\lambda}}} \in \mathbb{R}$ when $\left|N^\lambda - \rho^{\lambda}\right| \in \Theta(\!\sqrt{\lambda})\cup o(\!\sqrt{\lambda})$.
\end{theorem}
Section~\ref{sec:proof-of-theorem1} is devoted to the proof of Theorem~\ref{thm:ErlC-ASYM}.

%\blue{The well-known QED approximations~\cite{halfin1981heavy} fall within the fourth case outlined in Theorem~\ref{thm:ErlC-ASYM}, wherein the delay probability converges to a non-degenerate limit between $0$ and $1$. Our expanded understanding of the asymptotic properties of the Erlang-C formula under the other staffing rules extends its applicability to a wider range of scenarios in many-server queues.}

\section{Implications of Theorem~\ref{thm:ErlC-ASYM}}\label{sec:implication}

In understaffed systems (i.e., when $N^{\lambda} - \rho^{\lambda}$ is negative), the Erlang-C formula~\eqref{eq:ErlangC} loses its interpretation as the steady-state probability of delay in the $M$/$M$/$N$ queue but remains a well-defined mathematical expression that appears in the calculation of key performance measures (KPMs) in extensions of the $M$/$M$/$N$ queue (such as the $M$/$M$/$N$/$k$ and $M$/$M$/$N$+$M$ queues).
A central implication of Theorem~\ref{thm:ErlC-ASYM} for such systems is that it can be used to derive limiting approximations of these KPMs that help understand their dependence on the staffing rule.
Commonly used KPMs include the steady-state probability of delay~\cite{borst2004dimensioning,mandelbaum2009staffing}, probability of abandonment~\cite{garnett2002designing}, server utilization~\cite{zhong2023behavior}, and expected wait time~\cite{zeltyn2005call}, among others.
In this section, we demonstrate how Theorem~\ref{thm:ErlC-ASYM} can be leveraged to obtain the limiting value of the steady-state probability of delay in an $M$/$M$/$N$+$M$ queue (a setting in which square-root safety staffing is not always asymptotically optimal~\cite{mandelbaum2009staffing}) and briefly discuss the consequent insights regarding the appropriate choice of the staffing rule.

The $M$/$M$/$N$+$M$ queue (also known as the Erlang-A model, first introduced in~\cite{palm1943intensitatsschwankungen}) extends the $M$/$M$/$N$ queue by allowing customers waiting in queues to renege if they run out of patience before their service begins. The patience time of each customer is independent and identically distributed according to an exponential distribution with mean $1/\theta$. The steady-state probability of delay depends on the Erlang-C formula, as shown next in Lemma~\ref{lemma:prob-of-delay}.

\begin{lemma}\label{lemma:prob-of-delay}
	The steady-state probability of delay in the $M$/$M$/$N$+$M$ queue is given by
	\begin{align}\label{eq:prob-of-delay-ErlangA}
		P(N,\rho;\mu,\theta) = \left(1 + \frac{C(N,\rho)^{-1} - 1}{\mu(N-\rho) J(N,\rho;\mu,\theta)}\right)^{-1},
	\end{align}
	where
	\begin{align}\label{eq:J}
		J(N,\rho;\mu,\theta) := \int_{0}^{\infty} e^{\frac{\mu\rho}{\theta} \left(1 - e^{-\theta x}\right) - N\mu x}dx.
	\end{align}
	Furthermore, the quantity $\mu(N-\rho) J(N,\rho;\mu,\theta)$ also admits the following integral representation:
	\begin{align}\label{eq:J-quantity-integral}
		\mu(N-\rho) J(N,\rho;\mu,\theta) = 1 + \frac{\mu\rho}{\theta} \int_{0}^{1} e^{\frac{\mu\rho}{\theta}v} (1-v)^{\frac{N\mu}{\theta}-1} vdv.
	\end{align}
\end{lemma}
\begin{remark}
	Exact expressions for various performance measures in the $M$/$M$/$N$+$M$ and $M$/$M$/$N$+$G$ queues have also been derived in~\cite{garnett2002designing,zeltyn2005call}.
\end{remark}
The proof of Lemma~\ref{lemma:prob-of-delay} can be found in~\ref{appendix:proofs}.

Theorem~\ref{thm:ErlC-ASYM} can be leveraged to evaluate the limiting value of the steady-state probability of delay~\eqref{eq:prob-of-delay-ErlangA} under different staffing rules, as we show next in Proposition~\ref{prop:prob-of-delay}.

\begin{proposition}[\texorpdfstring{$\bm{M/M/N}$+$\bm{M}$}{M/M/N+M} Steady-State Delay Probability]\label{prop:prob-of-delay}
	\begin{align*}
		\lim_{\lambda\to\infty}& P(N^{\lambda},\rho^{\lambda};\mu,\theta) \\
		=&
		\begin{cases}
			1, & 0 < \rho^{\lambda} - N^{\lambda} \in \omega(\!\sqrt{\lambda}), \\
			\Big(1-\frac{\xi(z)}{\xi(-\hat{z})}\Big)^{-1} 
			\in 
			\left(\left(1+\!\sqrt{\frac{\theta}{\mu}}\right)^{-1}, 1\right), & 0 < \rho^{\lambda} - N^{\lambda} \in \Theta(\!\sqrt{\lambda}), \\
			\left(1+\!\sqrt{\frac{\theta}{\mu}}\right)^{-1}, & \left|N^{\lambda} - \rho^{\lambda}\right| \in o(\!\sqrt{\lambda}), \\
			\Big(1-\frac{\xi(z)}{\xi(-\hat{z})}\Big)^{-1} \in
			\left(0, \left(1+\!\sqrt{\frac{\theta}{\mu}}\right)^{-1}\right), & 0 < N^{\lambda} -\rho^{\lambda} \in \Theta(\!\sqrt{\lambda}), \\
			0, & 0 < N^{\lambda} -\rho^{\lambda} \in \omega(\!\sqrt{\lambda}),
		\end{cases}
	\end{align*}
	where $z = \lim_{\lambda\to\infty} \frac{\rho^{\lambda} - N^{\lambda}}{\sqrt{N^{\lambda}}}$ and $\hat{z} = z\!\sqrt{\frac{\mu}{\theta}}$.
\end{proposition}
The proof of Proposition~\ref{prop:prob-of-delay} is quite similar to that of Theorem~\ref{thm:ErlC-ASYM}, and is provided in~\ref{appendix:proofs}.

\begin{figure}[ht]
	\centering
	\includegraphics[scale=0.35]{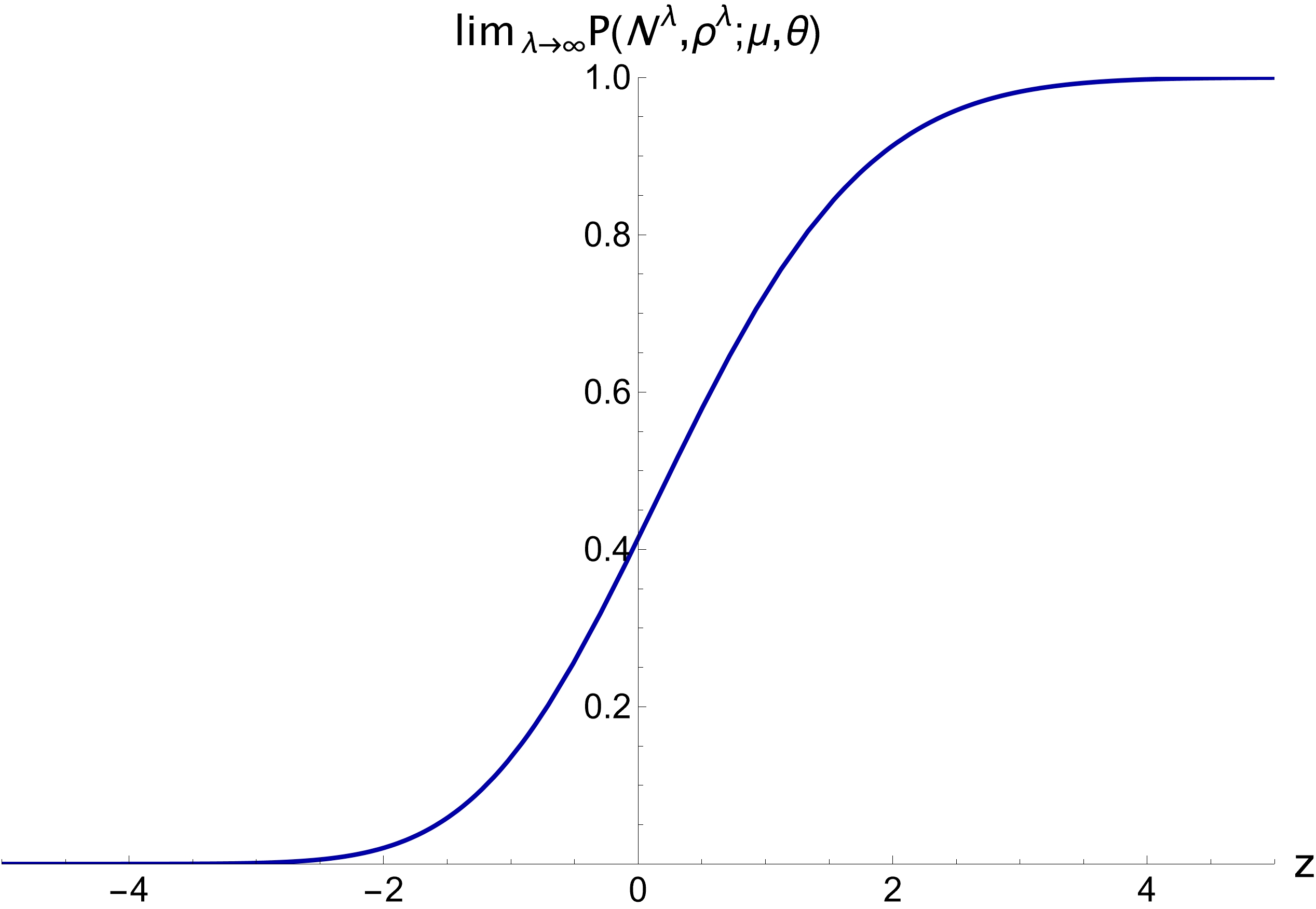}
	\caption{
		$\lim_{\lambda\to\infty} P(N^{\lambda},\rho^{\lambda};\mu,\theta)$ as a function of $z$ when $\mu=5$ and $\theta=10$.}
	\label{fig:limit-P}
\end{figure}

From Proposition~\ref{prop:prob-of-delay}, we see strict separation in the limiting steady-state probabilities of delay under different staffing rules in the critically loaded regime (where $\rho^{\lambda}/N^{\lambda}$ approaches~$1$ in the limit). Moreover, we numerically observe from Figure~\ref{fig:limit-P} that the limiting steady-state probability of delay appears to be concave in $|z|$ (for large enough $|z|$), meaning that changes in staffing exert the most significant impact on the probability of delay when the staffing rule is more balanced (i.e., $|z|$ is small), in a similar spirit as the law of diminishing marginal returns. This observation helps us to better understand and anticipate the trade-off between delay costs and staffing costs, which has managerial implications for the choice of optimal staffing rules.

\section{Proof of Theorem~\ref{thm:ErlC-ASYM}}\label{sec:proof-of-theorem1}

The proof leverages the following two auxiliary lemmas, whose proofs are deferred to the end of this section.% Section~\ref{sec:proof-of-lemmas}.
\begin{lemma}\label{lemma:integral-representation-C-inverse}
	The following are two equivalent integral representations of the reciprocal of $C(N,\rho)$:
	\begin{align}
		\mbox{(a)}\;\;& C(N,\rho)^{-1} = \rho\int_{0}^{\infty} e^{-\rho v}(1+v)^{N-1} v dv.\label{eq:C-inverse-new}\\
		\mbox{(b)}\;\;& C(N,\rho)^{-1} = 1+(N-\rho) \int_{-\infty}^{0} e^{-\rho(e^{-u}-1)-Nu} du.\label{eq:C-inverse}
	\end{align}
\end{lemma}

\begin{lemma}\label{lemma:asym-expansion-varphi}
	Let $\varphi(v;x) := e^{\frac{1}{2}v^2 - x\left(e^{-\frac{v}{\sqrt{x}}} + \frac{v}{\sqrt{x}} - 1\right)}$ for all $v \in \mathbb{R}$ and $x \in \mathbb{R}_{+}$. Then,
	\begin{align*}
		\varphi(v;x) = \sum_{i=0}^{\infty} a_i(v) \cdot x^{-\frac{i}{2}}\quad \forall v \in \mathbb{R}\ \forall x \in \mathbb{R}_{+},
	\end{align*}
	where $a_0(v)=1$ and $a_i(v)$ is a finite polynomial in $v$ for all $i\geq 1$.
\end{lemma}
We are now ready to prove Theorem~\ref{thm:ErlC-ASYM}.

\begin{itemize}[leftmargin=1.5em]
	\item \textbf{When}  $\left|N^\lambda - \rho^{\lambda}\right| \in \Theta(\!\sqrt{\lambda})\cup o(\!\sqrt{\lambda})$:
	Without loss of generality, let $N^\lambda = \rho^{\lambda} - z\sqrt{\rho^{\lambda}} + o(\!\sqrt{\rho^{\lambda}})$ for $z \in \mathbb{R}$. Note that $z = \lim_{\lambda\to\infty} \frac{\rho^{\lambda} - N^{\lambda}}{\sqrt{N^{\lambda}}}$.
	
	From Lemma~\ref{lemma:integral-representation-C-inverse}~(b),
	\begin{align*}
		C(N^{\lambda},\rho^{\lambda})^{-1} - 1 = (N^{\lambda}-\rho^{\lambda}) \int_{-\infty}^{0} e^{-\rho^{\lambda}(e^{-u}-1)-N^{\lambda} u} du.
	\end{align*}
	Let $v =\!\sqrt{\rho^{\lambda}} \ u$. The above can be equivalently written as
	\begin{align*}
		C(N^{\lambda}&,\rho^{\lambda})^{-1} - 1\\
		&= \left(\frac{N^{\lambda} - \rho^{\lambda}}{\sqrt{\rho^{\lambda}}}\right) \int_{-\infty}^{0} e^{-\left(\frac{N^{\lambda} - \rho^{\lambda}}{\sqrt{\rho^{\lambda}}}\right)v}  e^{-\rho^{\lambda}\left(e^{-\frac{v}{\sqrt{\rho^{\lambda}}}} + \frac{v}{\sqrt{\rho^{\lambda}}} -1\right)} dv\\
		&= \left(\frac{N^{\lambda} - \rho^{\lambda}}{\sqrt{\rho^{\lambda}}}\right) \int_{-\infty}^{0} e^{-\frac{1}{2}v^2-\left(\frac{N^{\lambda} - \rho^{\lambda}}{\sqrt{\rho^{\lambda}}}\right)v}  e^{\frac{1}{2}v^2-\rho^{\lambda}\left(e^{-\frac{v}{\sqrt{\rho^{\lambda}}}} + \frac{v}{\sqrt{\rho^{\lambda}}} -1\right)} dv\\
		&= \left(\frac{N^{\lambda} - \rho^{\lambda}}{\sqrt{\rho^{\lambda}}}\right)
		e^{\frac{1}{2}\left(\frac{N^{\lambda} - \rho^{\lambda}}{\sqrt{\rho^{\lambda}}}\right)^2}\int_{-\infty}^{0} e^{-\frac{1}{2}\left(v+\frac{N^{\lambda} - \rho^{\lambda}}{\sqrt{\rho^{\lambda}}}\right)^2}  \varphi(v;\rho^{\lambda}) dv.
	\end{align*}
	Substituting for $\varphi(v;\rho^{\lambda})$ from Lemma~\ref{lemma:asym-expansion-varphi} and taking the limit as $\lambda\to\infty$, we obtain:
	\begin{align*}
		\lim_{\lambda\to\infty}&C(N^{\lambda},\rho^{\lambda})^{-1} - 1 \\
		&= -z\ \phi(z)^{-1}\lim_{\lambda\to\infty}\sum_{i=0}^{\infty} \left(\rho^{\lambda}\right)^{-\frac{i}{2}} \int_{-\infty}^{0}
		e^{-\frac{1}{2}\left(v+\frac{N^{\lambda} - \rho^{\lambda}}{\sqrt{\rho^{\lambda}}}\right)^2} a_i(v) dv.
	\end{align*}
	Finally, noting that every integral within the sum is bounded (since $a_i(v)$ are all finite polynomials in $v$) and so, only the first term of the sum would survive in the limit, we obtain:
	\begin{align*}
		\lim_{\lambda\to\infty}C(N^{\lambda},\rho^{\lambda})^{-1} &= 1 - z\ \phi(z)^{-1}\lim_{\lambda\to\infty} \int_{-\infty}^{0}
		e^{-\frac{1}{2}\left(v+\frac{N^{\lambda} - \rho^{\lambda}}{\sqrt{\rho^{\lambda}}}\right)^2} dv \nonumber \\
		&= 1 - z\ \phi(z)^{-1}\int_{-\infty}^{0}
		e^{-\frac{1}{2}(v-z)^2} dv \nonumber \\
		&= 1 - z\ \phi(z)^{-1}\int_{z}^{\infty}
		e^{-\frac{1}{2}t^2} dt \nonumber \\
		&= 1 - z \frac{\Phi^c(z)}{\phi(z)} = 1-\xi(z),
	\end{align*}
	implying that
	\begin{align}\label{eq:limit-C-QED}
		\lim_{\lambda\to\infty} C(N^{\lambda}, \rho^{\lambda}) = \big(1-\xi(z)\big)^{-1}.
	\end{align}

	Furthermore, from Lemma~\ref{le:std-Normal} in~\ref{appendix:normal}:
	\begin{itemize}[leftmargin=2em]
		\item If $z>0$, which is the case when $0 < \rho^{\lambda} - N^{\lambda} \in \Theta(\!\sqrt{\lambda})$, then $\big(1-\xi(z)\big)^{-1} \in (1,\infty)$.
		\item If $z=0$, which is the case when $\left|N^{\lambda} - \rho^{\lambda}\right| \in o(\!\sqrt{\lambda})$, then $\big(1-\xi(z)\big)^{-1} = 1$.
		\item If $z<0$, which is the case when $0 < N^{\lambda} - \rho^{\lambda} \in \Theta(\!\sqrt{\lambda})$, then $\big(1-\xi(z)\big)^{-1} \in (0,1)$.
	\end{itemize}
	
	\item \textbf{When} $0 < \rho^{\lambda} -N^{\lambda} \in \omega(\!\sqrt{\lambda})$:
	Let $f(\lambda) = \rho^{\lambda} - N^{\lambda}$.
	Let $g(\lambda;z) = \rho^{\lambda} - z \sqrt{\rho^{\lambda}} + o(\!\sqrt{\rho^{\lambda}})$ for $z>0$, so that $0 < \rho^{\lambda} - g(\lambda;z) \in \Theta(\!\sqrt{\lambda})$.
	
	By definition of the $\omega$ and $\Theta$ notations, for any $z>0$, there exists $\Lambda(z)$ such that $f(\lambda) \geq \rho^{\lambda} - g(\lambda;z) > 0$ for all $\lambda \geq \Lambda(z)$.
	Since $C(N,\rho)$ is strictly decreasing in $N$~\cite[p.8]{WhittNotes}, it follows that, for all $\lambda \geq \Lambda(z)$,
	\begin{align*}
		C\left(N^{\lambda},\rho^{\lambda}\right) = C\left(\rho^{\lambda}-f(\lambda),\rho^{\lambda}\right) \geq C\left(g(\lambda;z),\rho^{\lambda}\right),
	\end{align*}
	implying that
	\begin{align*}
		\lim_{\lambda\to\infty} C\left(N^{\lambda},\rho^{\lambda}\right) \geq \lim_{\lambda\to\infty} C\left(g(\lambda;z),\rho^{\lambda}\right)
		= \big(1-\xi(z)\big)^{-1},
	\end{align*}
	from~\eqref{eq:limit-C-QED}. Since $z>0$ is chosen arbitrarily, we have
	\begin{align*}
		\lim_{\lambda\to\infty} C\left(N^{\lambda},\rho^{\lambda}\right) \geq \sup_{z > 0} \big(1-\xi(z)\big)^{-1} = \infty,
	\end{align*}
	from Lemma~\ref{le:std-Normal}~(b) in~\ref{appendix:normal}. Therefore,
	\begin{align*}
		\lim_{\lambda \rightarrow \infty} C\left(N^{\lambda},\rho^{\lambda}\right) = \infty.
	\end{align*}
	
	\item \textbf{When} $0 < N^{\lambda} - \rho^{\lambda} \in \omega(\!\sqrt{\lambda})$. Let $f(\lambda) = N^{\lambda} - \rho^{\lambda}$. Let $g(\lambda;z) = \rho^{\lambda} - z \sqrt{\rho^{\lambda}} + o(\!\sqrt{\rho^{\lambda}})$ for $z<0$, so that $0 < g(\lambda;z) - \rho^{\lambda} \in \Theta(\!\sqrt{\lambda})$.
	
	By definition of the $\omega$ and $\Theta$ notations, for any $z<0$, there exists $\Lambda(z)$ such that $f(\lambda) \geq g(\lambda;z) - \rho^{\lambda} > 0$ for all $\lambda \geq \Lambda(z)$.
	Since $C(N,\rho)$ is strictly decreasing in $N$~\cite[p.8]{WhittNotes}, it follows that, for all $\lambda \geq \Lambda(z)$,
	\begin{align*}
		C\left(N^{\lambda},\rho^{\lambda}\right) = C\left(\rho^{\lambda} + f(\lambda),\rho^{\lambda}\right) \leq C\left(g(\lambda;z),\rho^{\lambda}\right),
	\end{align*}
	implying that
	\begin{align*}
		\lim_{\lambda\to\infty} C\left(N^{\lambda},\rho^{\lambda}\right) \leq \lim_{\lambda\to\infty} C\left(g(\lambda;z),\rho^{\lambda}\right) = \big(1-\xi(z)\big)^{-1},
	\end{align*}
	from~\eqref{eq:limit-C-QED}. Since $z < 0$ is chosen arbitrarily, we have
	\begin{align*}
		\lim_{\lambda\to\infty} C\left(N^{\lambda},\rho^{\lambda}\right)
		\leq \inf_{z < 0} \big(1-\xi(z)\big)^{-1}
		= 0,
	\end{align*}
	from Lemma~\ref{le:std-Normal}~(b). Moreover, it follows from Lemma~\ref{lemma:integral-representation-C-inverse}(a) that $C(N,\rho)\geq 0$ for all $N$ and $\rho$. Therefore,
	\begin{align*}
		\lim_{\lambda \rightarrow \infty} C\left(N^{\lambda}, \rho^{\lambda}\right) = 0.
	\end{align*}\hfill$\blacksquare$
\end{itemize}

%\section{Proofs of Lemmas~\ref{lemma:integral-representation-C-inverse} and~\ref{lemma:asym-expansion-varphi}}\label{sec:proof-of-lemmas}

\subsection{Proof of Lemma~\ref{lemma:integral-representation-C-inverse}}
We begin by introducing a closely related performance measure, which is the steady-state probability of blocking in the $M$/$M$/$N$/$N$ queue (also known as the Erlang loss model), where a customer, upon arrival, is immediately lost if all the servers are busy. This quantity is given by the well-known Erlang-B formula~\citep[p.80]{Cooper72}:
\begin{align*}
	B\left(N,\rho\right) = \frac{\rho^N}{N!} \left(\sum_{i=0}^{N}\frac{\rho^i}{i!}\right)^{-1},
\end{align*}
whose reciprocal admits the following integral representation~\citep[Theorem~3]{jagerman1974some}:
\begin{align}\label{eq:Erlang-B-integral}
	B(N,\rho)^{-1} = \rho \int_{0}^{\infty} e^{-\rho v} (1+v)^{N} dv.
\end{align}
Also, the Erlang-B and Erlang-C formula are related~\citep[p.92]{Cooper72}:
\begin{align*}
	C\left(N,\rho\right) = \frac{N \cdot B(N,\rho)}{N - \rho \left(1 - B(N,\rho)\right)},
\end{align*}
which implies that
\begin{align}\label{eq:ErlangB-C}
	C(N,\rho)^{-1} = 1 + \left(1-\frac{\rho}{N}\right)\left(B(N,\rho)^{-1}-1\right).
\end{align}

We are now ready to prove Lemma~\ref{lemma:integral-representation-C-inverse}. First, the integral appearing on the right hand side of~\eqref{eq:C-inverse-new} can be evaluated as
\begin{align*}
	\rho\int_{0}^{\infty}& e^{-\rho v}(1+v)^{N-1} v dv\\
	&= \rho\int_{0}^{\infty} e^{-\rho v}(1+v)^{N-1} (1+v-1) dv\\
	&= \rho\int_{0}^{\infty} e^{-\rho v}(1+v)^{N} dv - \rho\int_{0}^{\infty} e^{-\rho v}(1+v)^{N-1} dv\\
	&= \rho\int_{0}^{\infty} e^{-\rho v}(1+v)^{N} dv - \frac{\rho}{N}\int_{0}^{\infty} e^{-\rho v}d\left((1+v)^{N}\right)\\
	&= \rho\int_{0}^{\infty} e^{-\rho v}(1+v)^{N} dv\\
	&\qquad - \frac{\rho}{N} \left(e^{-\rho v}(1+v)^{N}\Big\vert_{0}^{\infty}+\rho\int_{0}^{\infty} e^{-\rho v}(1+v)^{N} dv\right)\\
	&= B(N,\rho)^{-1} - \frac{\rho}{N}\left(-1 + B(N,\rho)^{-1}\right)\\
	&= 1 + \left(1-\frac{\rho}{N}\right)\left(B(N,\rho)^{-1}-1\right) = C(N,\rho)^{-1},
\end{align*}
from~\eqref{eq:ErlangB-C}, thereby establishing~\eqref{eq:C-inverse-new}.

Next, letting $e^{-u} = 1+v$, (\ref{eq:Erlang-B-integral}) can be equivalently written as
\begin{align}\label{eq:B-inverse}
	B(&N,\rho)^{-1} \nonumber\\
	&= \rho \int_{-\infty}^{0} e^{-\rho(e^{-u}-1)-Nu} \cdot e^{-u} du \nonumber\\
	&= \int_{-\infty}^{0} e^{-\rho(e^{-u}-1)-Nu} \left(\rho e^{-u} - N\right) du + N \int_{-\infty}^{0} e^{-\rho(e^{-u}-1)-Nu} du \nonumber\\
	%=& e^{-\rho(e^{-u}-1)-Nu} \Big|_{u=-\infty}^{0} + N \int_{-\infty}^{0} e^{-\rho(e^{-u}-1)-Nu} du \nonumber\\
	&= 1 + N \int_{-\infty}^{0} e^{-\rho(e^{-u}-1)-Nu} du.
\end{align}
Substituting for $B(N,\rho)^{-1}$ using (\ref{eq:B-inverse}) in~\eqref{eq:ErlangB-C} yields (\ref{eq:C-inverse}).\hfill$\blacksquare$

\subsection{Proof of Lemma~\ref{lemma:asym-expansion-varphi}}

We begin by noting that the Maclaurin series for the exponential function $e^t$ converges to the value of the function everywhere on its domain. We use it twice in the proof of Lemma~\ref{lemma:asym-expansion-varphi}; once for each of the two exponential functions in $\varphi(v;x)$. First, the exponent of $\varphi(v;x)$ can be evaluated as:
\begin{align*}
	\frac{1}{2}v^2 - x\,\Bigg(e^{-\frac{v}{\sqrt{x}}} &+ \frac{v}{\sqrt{x}} - 1\Bigg)\\
	&= \frac{1}{2}v^2 - x\left(\sum_{i=0}^{\infty}\frac{1}{i!}\left(-\frac{v}{\sqrt{x}}\right)^i + \frac{v}{\sqrt{x}} - 1\right)\\
	&= \frac{1}{2}v^2 - x\sum_{i=2}^{\infty}\frac{1}{i!}\left(-\frac{v}{\sqrt{x}}\right)^i\\
	&= - x\sum_{i=3}^{\infty}\frac{1}{i!}\left(-\frac{v}{\sqrt{x}}\right)^i\\
	&= \frac{v^3}{\sqrt{x}} \sum_{i=0}^{\infty} \frac{(-1)^i}{(i+3)!} \left(\frac{v}{\sqrt{x}}\right)^i.
\end{align*}
Using the Maclaurin series for $e^t$ once again, $\varphi(v;x)$ becomes:
\begin{align*}
	\varphi(&v;x) \\
	&= e^{\frac{v^3}{\sqrt{x}} \sum_{i=0}^{\infty} \frac{(-1)^i}{(i+3)!} \left(\frac{v}{\sqrt{x}}\right)^i} \\
	&= \sum_{j=0}^\infty\frac{1}{j!}\left(\frac{v^3}{\sqrt{x}} \sum_{i=0}^{\infty} \frac{(-1)^i}{(i+3)!} \left(\frac{v}{\sqrt{x}}\right)^i\right)^j \\
	&= 1 + \frac{v^3}{3!} x^{-\frac{1}{2}} - \left(\frac{v^4}{4!} - \frac{v^6}{2!(3!)^2}\right)x^{-1} + \left(\frac{v^5}{5!} - \frac{v^7}{3!4!} + \frac{v^9}{(3!)^4}\right) x^{-\frac{3}{2}} \\
	& \qquad\qquad\quad - \left(\frac{v^6}{6!} - \frac{9v^8}{2!4!5!} + \frac{v^{10}}{2!(3!)^2 4!} - \frac{v^{12}}{(3!)^4 4!}\right) x^{-2} + \ldots \\
	&= \sum_{i=0}^{\infty}a_i(v) \cdot x^{-\frac{i}{2}},
\end{align*} 
where $a_0(v)=1$ and $a_i(v)$ is a finite polynomial in $v$ (with degree $3i$) for all $i\geq 1$.\hfill$\blacksquare$

\appendix

\section{Properties of the Standard Normal Distribution}\label{appendix:normal}
In this section, we collect properties of the standard Normal distribution that are used in our proofs.
Recall that $\phi(x) = \frac{1}{\sqrt{2\pi}} e^{-\frac{1}{2}x^2}$ and $\Phi^c(x) = \int_{x}^{\infty} \phi(t) dt$ are the density function and the complementary cumulative distribution function of the standard normal distribution, respectively, and $\xi(x) := x\frac{\Phi^c(x)}{\phi(x)}$.

\begin{lemma}\label{le:std-Normal}
	The function $\xi: \mathbb{R} \mapsto \mathbb{R}$ satisfies the following:
	\begin{enumerate}
		\item [(a)] $\xi(x)$ is a strictly increasing function of $x$;
		\item [(b)] $\lim_{x \to -\infty} \xi(x) = -\infty$, $\xi(0) = 0$, and $\lim_{x \to \infty} \xi(x) = 1$.
	\end{enumerate}
\end{lemma}

\noindent\textit{Proof of Lemma~\ref{le:std-Normal}.}
We begin by noting that
\begin{align}\label{eq:Dphi}
	\phi'(x) = \frac{1}{\sqrt{2\pi}} e^{-\frac{x^2}{2}} (-x) = -x \phi(x).
\end{align}
Moreover, note that
\begin{align}\label{eq:power-phi}
	\lim_{x\to\infty} x^n \phi(x) =  \lim_{x\to\infty} \frac{1}{\sqrt{2\pi}} x^n e^{-\frac{x^2}{2}} = 0, \ \forall n \in \mathbb{N},
\end{align}
since exponential decay dominates polynomial growth. This further implies that
\begin{align}\label{eq:power-Phi-comple}
	& \lim_{x \to \infty} x^n \Phi^c(x) = \lim_{x \to \infty} \frac{\Phi^c(x)}{x^{-n}} \overset{(1)}{=} \lim_{x \to \infty} \frac{-\phi(x)}{-n x^{-n-1}} \nonumber\\
	=& \lim_{x \to \infty} \frac{x^{n+1}\phi(x)}{n} \overset{(2)}{=} 0, \ \forall n \in \mathbb{N},
\end{align}
where~(1) follows from L'H\^{o}pital's rule, and~(2) follows from~\eqref{eq:power-phi}. \\

\noindent {\bf Proof of (a): } Taking the derivative of $\xi(x)$ and using~\eqref{eq:Dphi} to substitute for $\phi'(x)$, we get:
\begin{align*}
	\xi'(x) =& \frac{\phi(x) \left(\Phi^c(x) - x \phi(x)\right) - x \Phi^c(x) \phi'(x)}{\phi(x)^2} \\
	=& \frac{\phi(x) \left(\Phi^c(x) - x \phi(x)\right) + x^2 \Phi^c(x) \phi(x)}{\phi(x)^2} \\
	=& \frac{(1+x^2) \Phi^c(x) - x \phi(x)}{\phi(x)}.
\end{align*}
Since $\phi(x) > 0$ for all $x$, it suffices to show that the numerator of the above display is strictly positive for all $x$. Define $\eta(x) := (1+x^2) \Phi^c(x) - x \phi(x)$.
Differentiating $\eta(x)$ once and using~\eqref{eq:Dphi} to substitute for $\phi'(x)$ yields:
\begin{align*}
	\eta'(x) =& 2x \Phi^c(x) - (1+x^2) \phi(x) - \phi(x) - x \phi'(x) \\
	=& 2x \Phi^c(x) - (1+x^2) \phi(x) - \phi(x) + x^2 \phi(x) \\
	=& 2 (x \Phi^c(x) - \phi(x)).
\end{align*}
Differentiating $\eta'(x)$ once and using~\eqref{eq:Dphi} to substitute for $\phi'(x)$ yields:
\begin{align*}
	\eta''(x) =& 2\left(\Phi^c(x) - x \phi(x) - \phi'(x)\right) \\
	=& 2\left(\Phi^c(x) - x \phi(x) + x \phi(x)\right) = 2\Phi^c(x) > 0,
\end{align*}
implying that $\eta'$ is strictly increasing. Since $\lim_{x \to \infty} \eta'(x) = 0$ from~\eqref{eq:power-Phi-comple}, it follows that $\eta'(x) < 0$ for all $x$, implying that $\eta$ is strictly decreasing. Finally, since $\lim_{x\to\infty} \eta(x) = 0$ from~\eqref{eq:power-phi} and~\eqref{eq:power-Phi-comple}, it follows that $\eta(x) > 0$ for all $x$. Hence, $\xi'(x) > 0$ for all $x$, implying that $\xi(x)$ is a strictly increasing function of $x$. \\

\noindent {\bf Proof of (b): } It is straightforward to see that
\begin{align*}
	\lim_{x \to -\infty} x\frac{\Phi^c(x)}{\phi(x)} = -\infty \quad \mbox{ and } \quad \lim_{x \to 0} x\frac{\Phi^c(x)}{\phi(x)} = 0.
\end{align*}
Finally, as $x \to \infty$, note that
\begin{align*}
	\lim_{x \to \infty} x\frac{\Phi^c(x)}{\phi(x)} &\overset{(3)}{=} \lim_{x \to \infty}\frac{\Phi^c(x) - x \phi(x)}{\phi'(x)} \overset{(4)}{=} \lim_{x \to \infty}\frac{\Phi^c(x) - x \phi(x)}{-x \phi(x)}\\
	&= 1 - \lim_{x \to \infty} \frac{\Phi^c(x)}{x \phi(x)} \overset{(5)}{=} 1 - \lim_{x \to \infty} \frac{-\phi(x)}{\phi(x) + x \phi'(x)}\\
	&\overset{(6)}{=} 1 - \lim_{x \to\infty} \frac{-\phi(x)}{\phi(x) -x^2 \phi(x)} = 1 - \lim_{x \to \infty} \frac{1}{x^2-1} = 1,
\end{align*}
where~(3) and~(5) follow from L'H\^{o}pital's rule, and~(4) and~(6) follow from~\eqref{eq:Dphi}.\hfill$\blacksquare$

\section{Proofs from Section~\ref{sec:implication}}\label{appendix:proofs}

\subsection{Proof of Lemma~\ref{lemma:prob-of-delay}}\label{sec:proof-of-lemma1}

The proof of~\eqref{eq:prob-of-delay-ErlangA} leverages the expressions for the steady-state probabilities of the $M$/$M$/$N$+$M$ queue developed in~\cite{baccelli1981queues}. From Equations (4.3) and (4.6) in~\cite{baccelli1981queues}, the steady-state probability of having $i$ customers in the system (either being served or waiting in queue), for $i\in\{0,1,2,\ldots,N-1\}$, is given by
\begin{align*}
	p_i(N,\rho;\mu,\theta) = \frac{\rho^i}{i!} \left(\sum_{i=0}^{N-1} \frac{\rho^i}{i!} + \frac{\rho^{N-1}}{(N-1)!} \lambda J(N,\rho;\mu,\theta)\right)^{-1}\\
	= \frac{\rho^i}{i!} \left(\sum_{i=0}^{N-1} \frac{\rho^i}{i!} + \frac{\rho^{N}}{N!} (N\mu) J(N,\rho;\mu,\theta)\right)^{-1}.
\end{align*}
Thus, the steady-state probability of delay is given by 
\begin{align*}
	P(N&,\rho;\mu,\theta) \\
	&= 1 - \sum_{i=0}^{N-1} p_i(N,\rho;\mu,\theta) \\
	&= 1 - \left(\sum_{i=0}^{N-1} \frac{\rho^i}{i!}\right) \left(\sum_{i=0}^{N-1} \frac{\rho^i}{i!} + \frac{\rho^{N}}{N!} (N\mu) J(N,\rho;\mu,\theta)\right)^{-1},
\end{align*}
which implies that
\begin{align*}
	P(N,\rho;\mu,\theta) &= \left(1 + \frac{\left(\sum_{i=0}^{N-1}\frac{\rho^i}{i!}\right)\Big/\frac{\rho^{N}}{N!}}{ (N\mu) J(N,\rho;\mu,\theta)}\right)^{-1} \\
	&= \left(1 + \frac{\left(1-\frac{\rho}{N}\right)\left(\sum_{i=0}^{N-1}\frac{\rho^i}{i!}\right)\Big/\frac{\rho^{N}}{N!}}{ \mu(N-\rho) J(N,\rho;\mu,\theta)}\right)^{-1} \\
	&= \left(1 + \frac{C(N,\rho)^{-1}-1}{\mu(N-\rho) J(N,\rho;\mu,\theta)}\right)^{-1},
\end{align*}
where the last step follows from~\eqref{eq:ErlangC}.

Next, for the proof of~\eqref{eq:J-quantity-integral}, we first use~\eqref{eq:J} to write:
\begin{align*}
	\mu(N-\rho)J(N,\rho;\mu,\theta) =\mu (N-\rho) \int_{0}^{\infty} e^{\frac{\mu\rho}{\theta} \left(1-e^{-\theta x}\right) - N\mu x} dx.
\end{align*}
Letting $u = -\theta x$, the above can be equivalently written as
\begin{align*}
	\mu(N-&\rho)J(N,\rho;\mu,\theta)\\
	&= \frac{\mu}{\theta} (N-\rho) \int_{-\infty}^{0} e^{\frac{\mu \rho}{\theta} \left(1-e^{u}\right) + \frac{N\mu}{\theta} u }du \\
	&= \int_{-\infty}^{0} e^{\frac{\mu \rho}{\theta} \left(1-e^{u}\right) + \frac{N\mu}{\theta} u} \left(\frac{N\mu}{\theta} - \frac{\mu\rho}{\theta} e^u +\frac{\mu\rho}{\theta}(e^{u}-1)\right) du \\
	&= \int_{-\infty}^{0} e^{\frac{\mu \rho}{\theta} \left(1-e^{u}\right) + \frac{N\mu}{\theta} u} \left(\frac{N\mu}{\theta} - \frac{\mu\rho}{\theta} e^u\right) du \\
	&\qquad\qquad\qquad+ \frac{\mu\rho}{\theta} \int_{-\infty}^{0} e^{\frac{\mu \rho}{\theta} \left(1-e^{u}\right) + \frac{N\mu}{\theta} u}   \left(e^u-1\right) du \\
	=& 1 +  \frac{\mu\rho}{\theta} \int_{-\infty}^{0} e^{\frac{\mu \rho}{\theta} \left(1-e^{u}\right) + \frac{N\mu}{\theta} u}   \left(e^u-1\right) du.
\end{align*}
Letting $v=1-e^u$, the above can be equivalently written as 
\begin{align*}
	\mu(N-\rho)J(N,\rho;\mu,\theta) &= 1 + \frac{\mu\rho}{\theta} \int_{1}^{0} e^{\frac{\mu\rho}{\theta}v + \frac{N\mu}{\theta} \log(1-v)} (-v) \left(\frac{dv}{1-v}\right) \\
	&= 1 + \frac{\mu\rho}{\theta} \int_{0}^{1} e^{\frac{\mu\rho}{\theta}v} (1-v)^{\frac{N\mu}{\theta}-1} vdv.\qquad\;\blacksquare
\end{align*}

\subsection{Proof of Proposition~\ref{prop:prob-of-delay}}\label{sec:proof-of-prop1}
The proof leverages the following auxiliary lemma, whose proof is deferred to the next subsection,~\ref{ssec:proof-off-applemma}.%\red{Proposition~\ref{prop:prob-of-delay} follows by applying Theorem~\ref{thm:ErlC-ASYM}, Lemma~\ref{lemma:prob-of-delay}, and the following auxiliary lemma.}

\begin{lemma}\label{lemma:J}
	If $N^{\lambda} = \rho^{\lambda} - z \sqrt{\rho^{\lambda}} + o(\!\sqrt{\rho^{\lambda}})$ for $z \in \mathbb{R}$, then $\lim_{\lambda\to\infty} \mu (N^{\lambda}-\rho^{\lambda}) J(N^{\lambda},\rho^{\lambda};\mu,\theta) = \xi(-\hat{z})$, where $\hat{z} = z\!\sqrt{\frac{\mu}{\theta}}$.
\end{lemma}

\begin{itemize}[leftmargin=1.5em]
	
	\item \textbf{When} $|N^{\lambda} - \rho^{\lambda}| \in \Theta(\!\sqrt{\lambda})\cup o(\!\sqrt{\lambda})$: Without loss of generality, let $N^{\lambda} = \rho^{\lambda} - z \sqrt{\rho^{\lambda}} + o(\!\sqrt{\rho^{\lambda}})$ for $z \in \mathbb{R}$. Note that $z = \lim_{\lambda\to\infty} \frac{\rho^{\lambda} - N^{\lambda}}{\sqrt{N^{\lambda}}}$. 
	From Lemma~\ref{lemma:prob-of-delay}, Theorem~\ref{thm:ErlC-ASYM}, and Lemma~\ref{lemma:J}, 
	\begin{align}\label{eq:prob-of-delay-QED}
		\lim_{\lambda\to\infty}P(N^{\lambda},\rho^{\lambda};\mu,\theta)
		&= \lim_{\lambda\to\infty} \left(1 + \frac{C(N^{\lambda},\rho^{\lambda})^{-1} - 1}{\mu(N^{\lambda}-\rho^{\lambda}) J(N^{\lambda},\rho^{\lambda};\mu,\theta)}\right)^{-1} \nonumber\\
		&= \left(1 - \frac{\xi(z)}{\xi(-\hat{z})}\right)^{-1}.
	\end{align}
	Furthermore, from Lemma~\ref{le:std-Normal} in~\ref{appendix:normal}, $\frac{\xi(z)}{\xi(-\hat{z})}$ is strictly increasing in $z$ from $-\infty$ (when $z\to-\infty$) to $-\!\sqrt{\frac{\theta}{\mu}}$ (at $z=0$) to $0$ (when $z\to\infty$), implying the following:
	\begin{itemize}
		\item If $z>0$, which is the case when $0 < \rho^{\lambda} - N^{\lambda} \in \Theta(\!\sqrt{\lambda})$, then $\left(1 - \frac{\xi(z)}{\xi(-\hat{z})}\right)^{-1} \in \left(\left(1+\sqrt{\frac{\theta}{\mu}}\right)^{-1},1\right)$. 
		\item If $z=0$, which is the case when $|N^{\lambda} - \rho^{\lambda}| \in o(\!\sqrt{\lambda})$, then $\left(1 - \frac{\xi(z)}{\xi(-\hat{z})}\right)^{-1} = 1$.
		\item If $z<0$, which is the case when $0 < N^{\lambda} - \rho^{\lambda} \in \Theta(\!\sqrt{\lambda})$, then $\left(1 - \frac{\xi(z)}{\xi(-\hat{z})}\right)^{-1} \in \left(0,\left(1+\sqrt{\frac{\theta}{\mu}}\right)^{-1}\right)$. 
	\end{itemize}

	\item \textbf{When} $0 < \rho^{\lambda} - N^{\lambda} \in \omega(\!\sqrt{\lambda})$:
	Let $f(\lambda) = \rho^{\lambda} - N^{\lambda}$.
	Let $g(\lambda;z) = \rho^{\lambda} - z \sqrt{\rho^{\lambda}} + o(\!\sqrt{\rho^{\lambda}})$ for $z>0$, so that $0 < \rho^{\lambda} - g(\lambda;z) \in \Theta(\!\sqrt{\lambda})$.
	
	By definition of the $\omega$ and $\Theta$ notations, for any $z>0$, there exists $\Lambda(z)$ such that $f(\lambda) \geq \rho^{\lambda} - g(\lambda;z) > 0$, for all $\lambda \geq \Lambda(z)$. Since $P(N,\rho;\mu,\theta)$, given by Lemma~\ref{lemma:prob-of-delay}, is strictly decreasing in $N$ (recalling that $C(N,\rho)$ is strictly decreasing in $N$~\cite[p.8]{WhittNotes} and the quantity $\mu(N-\rho)J(N,\rho;\mu,\theta)$ is also strictly decreasing in $N$ from (\ref{eq:J-quantity-integral})), it follows that, for all $\lambda \geq \Lambda(z)$,
	\begin{align*}
		P\left(N^{\lambda},\rho^{\lambda};\mu,\theta\right) &= P\left(\rho^{\lambda}-f(\lambda),\rho^{\lambda};\mu,\theta\right) \\
		&\geq P\left(g(\lambda;z),\rho^{\lambda};\mu,\theta\right),
	\end{align*}
	implying that
	\begin{align*}
		\lim_{\lambda\to\infty} P\left(N^{\lambda},\rho^{\lambda};\mu,\theta\right)
		&\geq \lim_{\lambda\to\infty} P\left(g(\lambda;z),\rho^{\lambda};\mu,\theta\right) \\
		&= \left(1 - \frac{\xi(z)}{\xi(-\hat{z})}\right)^{-1},
	\end{align*}
	from (\ref{eq:prob-of-delay-QED}). Since $z>0$ is chosen arbitrarily, we have 
	\begin{align*}
		\lim_{\lambda\to\infty} P\left(N^{\lambda},\rho^{\lambda};\mu,\theta\right)
		\geq \sup_{z>0} \left(1 - \frac{\xi(z)}{\xi(-\hat{z})}\right)^{-1} = 1, 
	\end{align*}
	from the preceding discussion. Therefore, 
	\begin{align*}
		\lim_{\lambda\to\infty}P(N^{\lambda},\rho^{\lambda};\mu,\theta) = 1.
	\end{align*}

	\item \textbf{When} $0 < N^{\lambda} - \rho^{\lambda} \in \omega(\!\sqrt{\lambda})$. Let $f(\lambda) = N^{\lambda} - \rho^{\lambda}$. Let $g(\lambda;z) = \rho^{\lambda} - z \sqrt{\rho^{\lambda}} + o(\!\sqrt{\rho^{\lambda}})$ for $z<0$, so that $0 < g(\lambda;z) - \rho^{\lambda} \in \Theta(\!\sqrt{\lambda})$.
	
	By definition of the $\omega$ and $\Theta$ notations, for any $z<0$, there exists $\Lambda(z)$ such that $f(\lambda) \geq g(\lambda;z) - \rho^{\lambda} > 0$, for all $\lambda \geq \Lambda(z)$. Since $P(N,\rho;\mu,\theta)$, given by Lemma~\ref{lemma:prob-of-delay}, is strictly decreasing in $N$ (recalling that $C(N,\rho)$ is strictly decreasing in $N$~\cite[p.8]{WhittNotes} and the quantity $\mu(N-\rho)J(N,\rho;\mu,\theta)$ is also strictly decreasing in $N$ from (\ref{eq:J-quantity-integral})), it follows that, for all $\lambda \geq \Lambda(z)$,
	\begin{align*}
		P\left(N^{\lambda},\rho^{\lambda};\mu,\theta\right) &= P\left(\rho^{\lambda} + f(\lambda),\rho^{\lambda};\mu,\theta\right) \\
		&\leq P\left(g(\lambda;z),\rho^{\lambda};\mu,\theta\right),
	\end{align*}
	implying that
	\begin{align*}
		\lim_{\lambda\to\infty} P\left(N^{\lambda},\rho^{\lambda};\mu,\theta\right)
		&\leq \lim_{\lambda\to\infty} P\left(g(\lambda;z),\rho^{\lambda};\mu,\theta\right) \\
		&= \left(1 - \frac{\xi(z)}{\xi(-\hat{z})}\right)^{-1},
	\end{align*}
	from (\ref{eq:prob-of-delay-QED}). Since $z<0$ is chosen arbitrarily, we have 
	\begin{align*}
		\lim_{\lambda\to\infty} P\left(N^{\lambda},\rho^{\lambda};\mu,\theta\right)
		\leq \inf_{z<0} \left(1 - \frac{\xi(z)}{\xi(-\hat{z})}\right)^{-1} = 0, 
	\end{align*}
	from the preceding discussion. Therefore, 
	\begin{align*}
		\lim_{\lambda\to\infty}P(N^{\lambda},\rho^{\lambda};\mu,\theta) = 0.
	\end{align*}\hfill$\blacksquare$
\end{itemize}

\subsection{Proof of Lemma~\ref{lemma:J}}\label{ssec:proof-off-applemma}
From (\ref{eq:J}), 
\begin{align*}
	\mu(N^{\lambda}-\rho^{\lambda})J(N^{\lambda},\rho^{\lambda};\mu,\theta) =\mu (N^{\lambda}-\rho^{\lambda}) \int_{0}^{\infty} e^{\frac{\mu\rho^{\lambda}}{\theta} \left(1-e^{-\theta x}\right) - N^{\lambda} \mu x} dx.
\end{align*}
Letting $u = \theta x$, the above can be equivalently written as
\begin{align*}
	\mu(N^{\lambda}-\rho^{\lambda})J(N^{\lambda},\rho^{\lambda};\mu,\theta) 
	= \frac{\mu}{\theta} (N^{\lambda}-\rho^{\lambda}) \int_{0}^{\infty} e^{\frac{\mu \rho^{\lambda}}{\theta} \left(1-e^{-u}\right) - \frac{N^{\lambda} \mu}{\theta} u }du.
\end{align*}
Let $\hat{N}^{\lambda}=\frac{\mu}{\theta}N^{\lambda}$ and $\hat{\rho}^{\lambda}=\frac{\mu}{\theta}\rho^{\lambda}$ for all $\lambda>0$, and $\hat{z}=z\!\sqrt{\frac{\mu}{\theta}}$. Then, $N^{\lambda} = \rho^{\lambda} - z \sqrt{\rho^{\lambda}} + o(\!\sqrt{\rho^{\lambda}})$ is equivalent to $\hat{N}^{\lambda} = \hat{\rho}^{\lambda} - \hat{z} \sqrt{\hat{\rho}^{\lambda}} + o(\!\sqrt{\hat{\rho}^{\lambda}})$ and the above can be equivalently written as 
\begin{align*}
	\mu(N^{\lambda}-\rho^{\lambda})J(N^{\lambda},\rho^{\lambda};\mu,\theta) 
	= (\hat{N}^{\lambda}-\hat{\rho}^{\lambda}) \int_{0}^{\infty} e^{-\hat{\rho}^{\lambda} \left(e^{-u}-1\right) - \hat{N}^{\lambda}u} du.
\end{align*}
Let $v =\!\sqrt{\hat{\rho}^{\lambda}} \ u$. The above can be equivalently written as
\begin{align*}
	\mu(N^{\lambda}&-\rho^{\lambda})J(N^{\lambda},\rho^{\lambda};\mu,\theta)\\
	&= \left(\frac{\hat{N}^{\lambda} - \hat{\rho}^{\lambda}}{\sqrt{\hat{\rho}^{\lambda}}}\right) \int_{0}^{\infty} e^{-\left(\frac{\hat{N}^{\lambda} - \hat{\rho}^{\lambda}}{\sqrt{\hat{\rho}^{\lambda}}}\right)v}  e^{-\hat{\rho}^{\lambda}\left(e^{-\frac{v}{\sqrt{\hat{\rho}^{\lambda}}}} + \frac{v}{\sqrt{\hat{\rho}^{\lambda}}} -1\right)} dv\\
	%&= \left(\frac{N^{\lambda} - 	\hat{\rho}^{\lambda}}{\sqrt{\hat{\rho}^{\lambda}}}\right) \int_{0}^{\infty} e^{-\frac{1}{2}v^2-\left(\frac{\hat{N}^{\lambda} - \hat{\rho}^{\lambda}}{\sqrt{\hat{\rho}^{\lambda}}}\right)v}  e^{\frac{1}{2}v^2-\hat{\rho}^{\lambda}\left(e^{-\frac{v}{\sqrt{\hat{\rho}^{\lambda}}}} + \frac{v}{\sqrt{\hat{\rho}^{\lambda}}} -1\right)} dv\\
	&= \left(\frac{\hat{N}^{\lambda} - \hat{\rho}^{\lambda}}{\sqrt{\hat{\rho}^{\lambda}}}\right)
	e^{\frac{1}{2}\left(\frac{\hat{N}^{\lambda} - \hat{\rho}^{\lambda}}{\sqrt{\hat{\rho}^{\lambda}}}\right)^2}\int_{0}^{\infty} e^{-\frac{1}{2}\left(v+\frac{\hat{N}^{\lambda} - \hat{\rho}^{\lambda}}{\sqrt{\hat{\rho}^{\lambda}}}\right)^2}  \varphi(v;\hat{\rho}^{\lambda}) dv.
\end{align*}
Substituting for $\varphi(v;\hat{\rho}^{\lambda})$ from Lemma~\ref{lemma:asym-expansion-varphi} and taking the limit as $\lambda\to\infty$, we obtain:
\begin{align*}
	\lim_{\lambda\to\infty}&\mu(N^{\lambda}-\rho^{\lambda})J(N^{\lambda},\rho^{\lambda};\mu,\theta) \\
	&= -\hat{z}\ \phi(-\hat{z})^{-1}\lim_{\lambda\to\infty}\sum_{i=0}^{\infty} \left(\hat{\rho}^{\lambda}\right)^{-\frac{i}{2}} \int_{0}^{\infty}
	e^{-\frac{1}{2}\left(v+\frac{\hat{N}^{\lambda} - \hat{\rho}^{\lambda}}{\sqrt{\hat{\rho}^{\lambda}}}\right)^2} a_i(v) dv,
\end{align*}
recalling that $\hat{z} = \lim_{\lambda\to\infty} \frac{\hat{\rho}^{\lambda} - \hat{N}^{\lambda}}{\sqrt{\hat{\rho}^{\lambda}}}$. Finally, noting that every integral within the sum is bounded (since $a_i(v)$ are all finite polynomials in $v$) and so, only the first term of the sum would survive in the limit, we obtain:
\begin{align*}
	\lim_{\lambda\to\infty}\mu(N^{\lambda}&-\rho^{\lambda})J(N^{\lambda},\rho^{\lambda};\mu,\theta)\\
	&= - \hat{z}\ \phi(-\hat{z})^{-1}\lim_{\lambda\to\infty} \int_{0}^{\infty}
	e^{-\frac{1}{2}\left(v+\frac{\hat{N}^{\lambda} - \hat{\rho}^{\lambda}}{\sqrt{\hat{\rho}^{\lambda}}}\right)^2} dv \\
	&= - \hat{z}\ \phi(-\hat{z})^{-1}\int_{0}^{\infty}
	e^{-\frac{1}{2}(v-\hat{z})^2} dv \\
	&= - \hat{z}\ \phi(-\hat{z})^{-1}\int_{-\hat{z}}^{\infty}
	e^{-\frac{1}{2}t^2} dt \\
	&= - \hat{z} \frac{\Phi^c(-\hat{z})}{\phi(-\hat{z})} = \xi(-\hat{z}).\qquad\qquad\qquad\qquad\qquad\quad\;\;\,\blacksquare
\end{align*}

\section*{Acknowledgments}
Ragavendran Gopalakrishnan is supported in part by funding from the Social Sciences and Humanities Research Council of Canada.

\bibliographystyle{ormsv080}
\bibliography{ErlangAsymProp}

\end{document}